\documentclass{amsart}
\usepackage{graphicx}
\usepackage{amsmath}
\usepackage{amsfonts}
\usepackage{amssymb}
\usepackage{tikz-cd} 	
\usepackage{parskip}
\usepackage{amsthm}
\usepackage{tikz}
\usetikzlibrary{calc}
\usepackage{times}
\usepackage{lineno}
\usepackage[super]{nth}
\graphicspath{{./figures/}}
\usepackage{etoolbox}

\title{The Buckyball Spectrum}
\author{G. Henderson-Walshe}
\date{June 2026}

\DeclareMathOperator{\id}{\mathbf{1}}

\DeclareMathOperator{\End}{End}

\DeclareMathOperator{\Mat}{Mat}

\DeclareMathOperator{\tr}{tr}

\DeclareMathOperator{\ZZ}{\mathbf{Z}}
\DeclareMathOperator{\NN}{\mathbf{N}}

\DeclareMathOperator{\CC}{\mathbf{C}}

\DeclareMathOperator{\cH}{\mathcal{H}}

\theoremstyle{definition}

\newcommand{\TitlePageQRCodeFile}{QR_code.jpg}
\newcommand{\TitlePageQRCodeSize}{2.2cm}
\newcommand{\TitlePageQRCodeXOffset}{2.5cm}
\newcommand{\TitlePageQRCodeYOffset}{2.5cm}

\let\oldmaketitle\maketitle

\renewcommand{\maketitle}{%
  \begin{tikzpicture}[remember picture, overlay]
    \node[
      anchor=north west,
      inner sep=0pt,
      outer sep=0pt
    ] at (
      $ (current page.north west)
      +(\TitlePageQRCodeXOffset,-\TitlePageQRCodeYOffset) $
    ) {%
      \includegraphics[width=\TitlePageQRCodeSize]{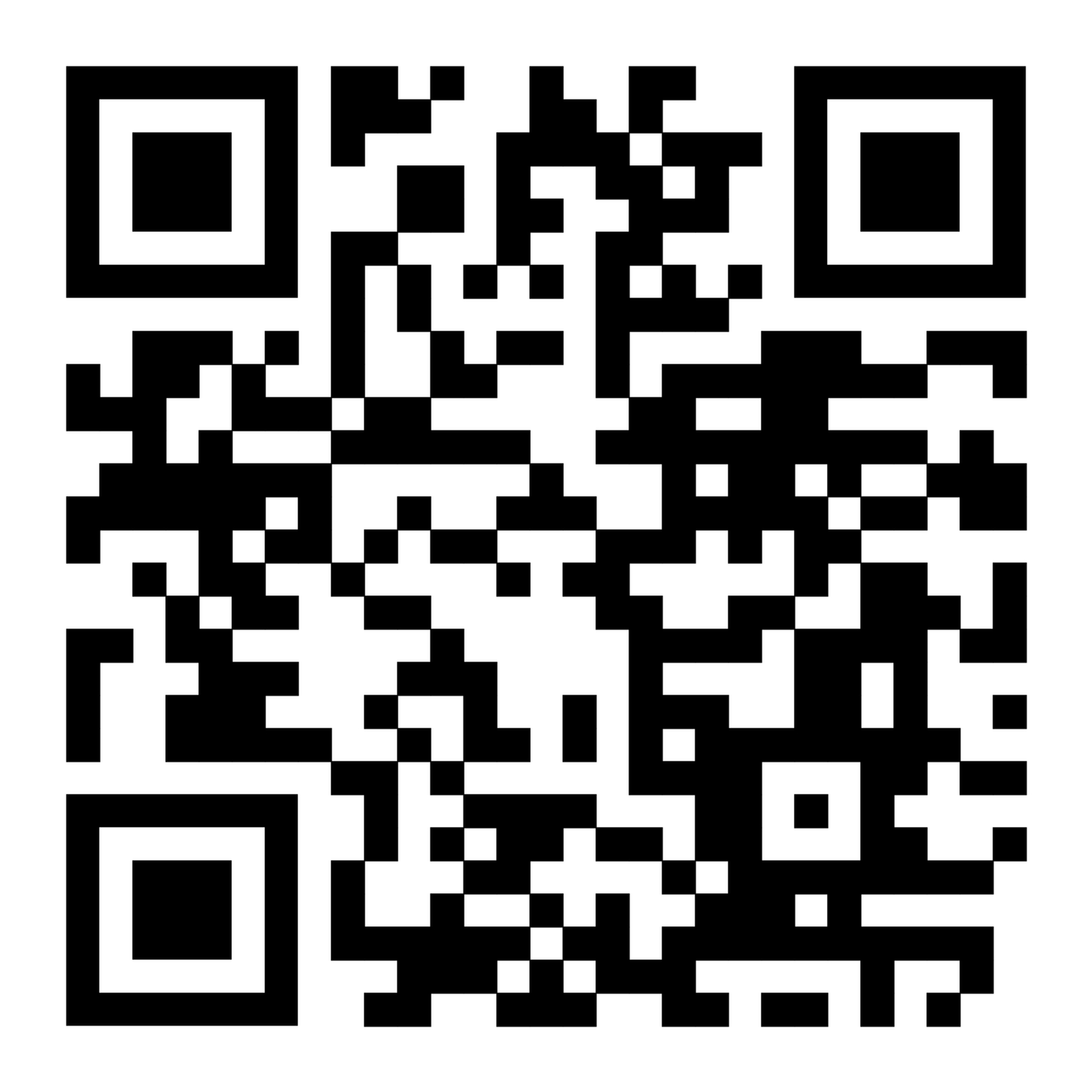}%
    };
  \end{tikzpicture}%
  \oldmaketitle
}

\begin{document}

\maketitle

In this article, we use representation theory to calculate the $\pi$-electron energy levels of Buckminsterfullerene and other symmetrical conjugated systems according to the H\"uckel approximation. We discuss the general case of a free transitive group action on a graph by graph automorphisms. We use representation theory to calculate the eigenvalues of the adjacency matrices, without needing to compute the matrix at all. The adjacency matrix is shown to be related to an element of the regular representation whose image under the irreducible characters is easy to calculate. The eigenvalues are then calculated from the characters using Newton's formula. Other examples, including cyclic hydrocarbons and infinite graphene sheets, are also discussed.

The article is accompanied by two supplementary Jupyter notebooks with code written in Sage, linked with the QR code in the top left. \textbf{Notebook 1} is used as a check, and computes the actual matrix by first computing matrices for the regular representation of the group. \textbf{Notebook 2} uses the characters of the group along with Newton's formula to find the eigenvalues without computing any matrices.

\section{Introduction: Hückel Theory}

In quantum chemistry, a \emph{state} is an element of the Hilbert space $\mathcal{H}$ of square-integrable $\CC$-valued functions on a region $\Omega$ (which we may rescale to have unit volume, $|\Omega|=1$), together with the inner product, $\langle \psi, \phi \rangle = \int_\Omega \psi^* \phi$, (where $*$ denotes complex conjugation). 

An \emph{observable} is an operator $T$ on $\cH$ such that $\langle \psi, T \psi \rangle$ is a measurable quantity. The function $\psi^* T \psi$ is the density function of the observable for that state. An electronic state is an element $\psi \in \mathcal{H}$ such that the density of the identity operator for $\psi$, $\psi^* \psi$, is the probability density function of an electron in the system, so $\langle \psi, \psi \rangle = 1$. The \emph{Hamiltonian} $\check{H}$ is the operator for which $\langle \psi, \check{H} \psi \rangle$ is the total energy of the state $\psi$. As total energy is conserved, the energy density should be proportional to the density of the state. In equations,

\begin{equation}
	\check{H} \psi = E \psi,
\end{equation}
where $E$ is the total energy of the state. You may recognise the form of this equation as telling us that the observable energy is an eigenvalue of the operator $\check{H}$. 

Now, to explain the type of quantum mechanical system which is the object of this article, we briefly overview some chemistry. The first two types of electron orbitals at an atom are $s$—which are spherical—and $p$—which are ``dumbbell''-shaped, orthogonal to the plane of the molecule (or rather the local approximation of a plane, as may be the case). Now, the orbitals are arranged into shells; shells which are filled have lower energy, and therefore are not available for bonding. Conversely, the valence shell may need to be filled by bonding. Atoms form $\sigma$-bonds when they have orbitals overlapping along the axis between the two atoms, while $\pi$-bonds are formed by orbitals which overlap parallel to the axis between the atoms, and therefore may only be formed by $p$ electrons. If there are not enough $p$ electrons to form the necessary $\sigma$-bonds, $s$ electrons must be brought up to a higher energy level; $s$ and $p$ orbitals are combined to form \emph{hybrid} orbitals.

For example, a carbon atom in its ground state has four $p$ electrons and two $s$ electrons. If it has four $\sigma$-bonds, instead of filling  its $\sigma$-bonds with two $p$ orbitals and two $s$ orbitals, it will have four ``$sp^3$ hybridised'' orbitals. If the carbon atom has three $\sigma$-bonds and one $\pi$-bond, corresponding to one double bond and three single bonds, the $\pi$-bond will be formed by a $p$ electron and the $\sigma$-bonds will be filled by $sp^2$ hybridised electrons.

When two pairs of $\pi$-bonded atoms are linked by a $sp$ or $sp^2$ hybridised $\sigma$-bond, the $\pi$-orbitals join together to form a \emph{conjugated system}, where all of the $\pi$-bonded electrons are \emph{de-localised}. Examples include benzene, naphthalene, graphene, butadiene, and of course, Buckminsterfullerene. \textbf{H\"uckel theory} is a method of approximating the configurations of the $\pi$-bonded electrons in a conjugated system. We will focus on molecules composed of only carbon and hydrogen atoms, where all of the carbons form one large conjugated system. Thus the number of $\pi$-bonded electrons in the system is equal to the number of carbons, which we will call $n$.

The idea is to approximate $\cH$ by a finite vector space corresponding to a basis of orbitals around each carbon in an organic molecule, and to approximate the Hamiltonian $\check{H}$ by a finite matrix
\begin{equation}
H = \alpha I + \beta M,
\end{equation}
where $\alpha$ and $\beta$ are real constants determined empirically, $I$ is the identity matrix, and $M$ is the adjacency matrix of the carbons in the molecule, considered as a graph. The adjacency matrix of a simple graph is just the matrix 
\begin{equation} \label{eq:adjacency matrix def}
(M)_{ij} = \begin{cases}
1 \text{ if the vertices $v_i, v_j$ share an edge in the graph} \\
0 \text{ if the vertices $v_i, v_j$ do not share an edge in the graph}
\end{cases}.
\end{equation}
For the Hückel theory version, substitute the words `graph' for `molecule', `vertex' for `carbon', and `edge' for `bond' in the previous sentence. For example, the molecule benzene has the adjacency matrix

\begin{equation} \label{eq:benzene huckel matrix}
 \begin{bmatrix}
0 & 1 & 0 & 0 & 0 & 1 \\
1 & 0 & 1 & 0 & 0 & 0 \\
0 & 1 & 0 & 1 & 0 & 0 \\
0 & 0 & 1 & 0 & 1 & 0 \\
0 & 0 & 0 & 1 & 0 & 1 \\
1 & 0 & 0 & 0 & 1 & 0 
\end{bmatrix} 
\end{equation}

\begin{figure} 
    \centering
    \includegraphics[width=0.5\linewidth]{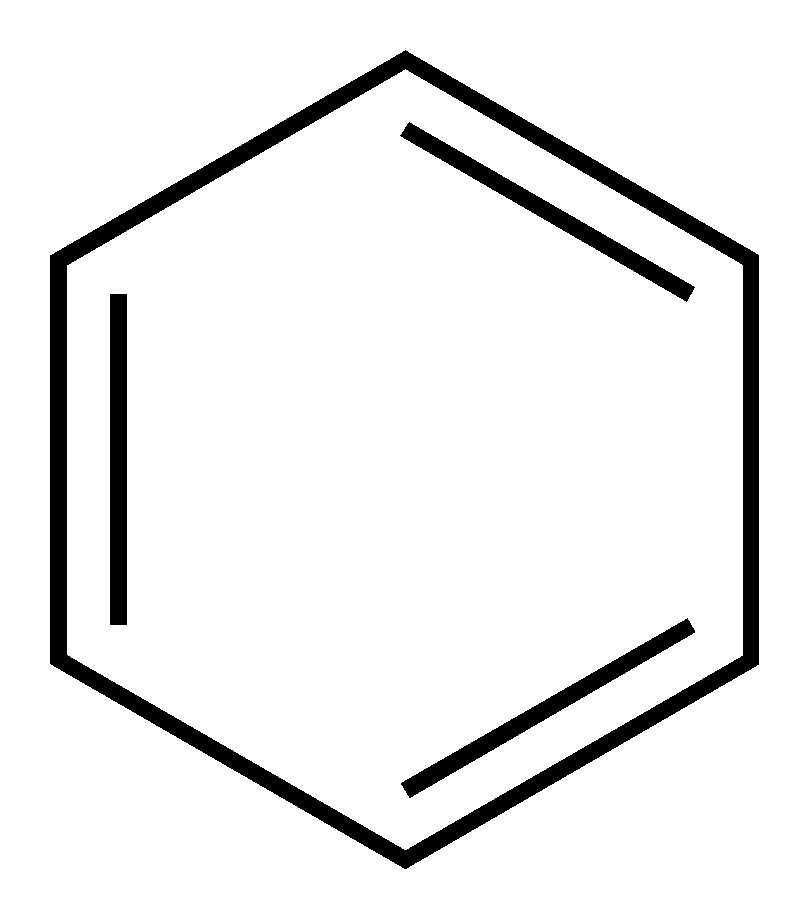}
    \caption{Skeleton of benzene. Adapted from \cite{benzene_image}.}
    \label{fig:benzene}
\end{figure}

You may notice that benzene has double as well as single bonds, and indeed we may wish to make the model more accurate by splitting the adjacency matrix into single and double bond terms $M_1$ and $M_2$,
\begin{equation} \label{eq:Huckel_matrix}
	H = \alpha I + \beta M_1 + \gamma M_2.
\end{equation} 

The ground state of the $\pi$-electrons is the assignment of the $n$ $\pi$-bonding electrons to $n$ pure states such that the total energy of the electrons is minimised. The Pauli-Exclusion Principle says that there can be at most two electrons in each state, so the multiset of possible energy levels is the multiset of energy levels counted twice. Therefore the ground state energy is the sum of the lowest $n$ of the multiset of eigenvalues, where each eigenvalue is counted twice. Lower ground state energy corresponds to a more stable molecule. For the purpose of comparing relative stability, we assume that $\alpha = 0$ and $\beta = \gamma = -1$ in Equation (\ref{eq:Huckel_matrix}).

\section{Free Transitive Group Actions on a Graph}

The general problem we are going to solve in this article is to find the eigenvalues of the adjacency matrix of a graph $X = (V,E)$ given a \textbf{free}, \textbf{transitive} action of a finite group $G$ by \textbf{graph automorphisms} of $X$. Let me explain what each of those terms means. We will abbreviate the adjectives `free and transitive by graph automorphisms' to `f.t.a.'.

An action of a group $G$ on a set $X$ is said to be \emph{transitive} if for each $x,x' \in X$ there exists a $g \in G$ such that $g x = x'$ and \emph{free} if such a $g$ is always unique. Therefore if $G$ acts freely and transitively on $X$ then for each $x \in X$ we get a bijection $G \leftrightarrow X$. An \emph{automorphism} of a graph $X=(V,E)$ is a bijection $f: V \to V$ such that for each edge $(v, v') \in E$, $(f^{-1}(v),f^{-1}(v')) \in E$ and $ (f(v),f(v')) \in E$. 

An example of an f.t.a. action is the action of the cyclic group $\ZZ/6\!\ZZ$ on the (graph of the) molecule benzene. It is classical that every finite subgroup of $SO(3)$ is a subgroup of a cyclic group or the (rotational) symmetry group of a regular tetrahedron, octahedron, or icosahedron \cite{zimmermann2011finitegroupsactingspheres}. It turns out that the symmetry group of the icosahedron is isomorphic to the alternating group $A_5$, which should be believable given that $|A_5|=5!/2=60$ and the symmetry group of the icosahedron has order $20 \times 3=60$ by the orbit-stabiliser theorem applied to a vertex. The Buckminsterfullerene molecule also has icosahedral symmetry, but each vertex's only stabiliser is the identity, so the action of the icosahedral group, isomorphic to $A_5$, is f.t.a.

\begin{figure}
    \centering
    \includegraphics[width=0.5\linewidth]{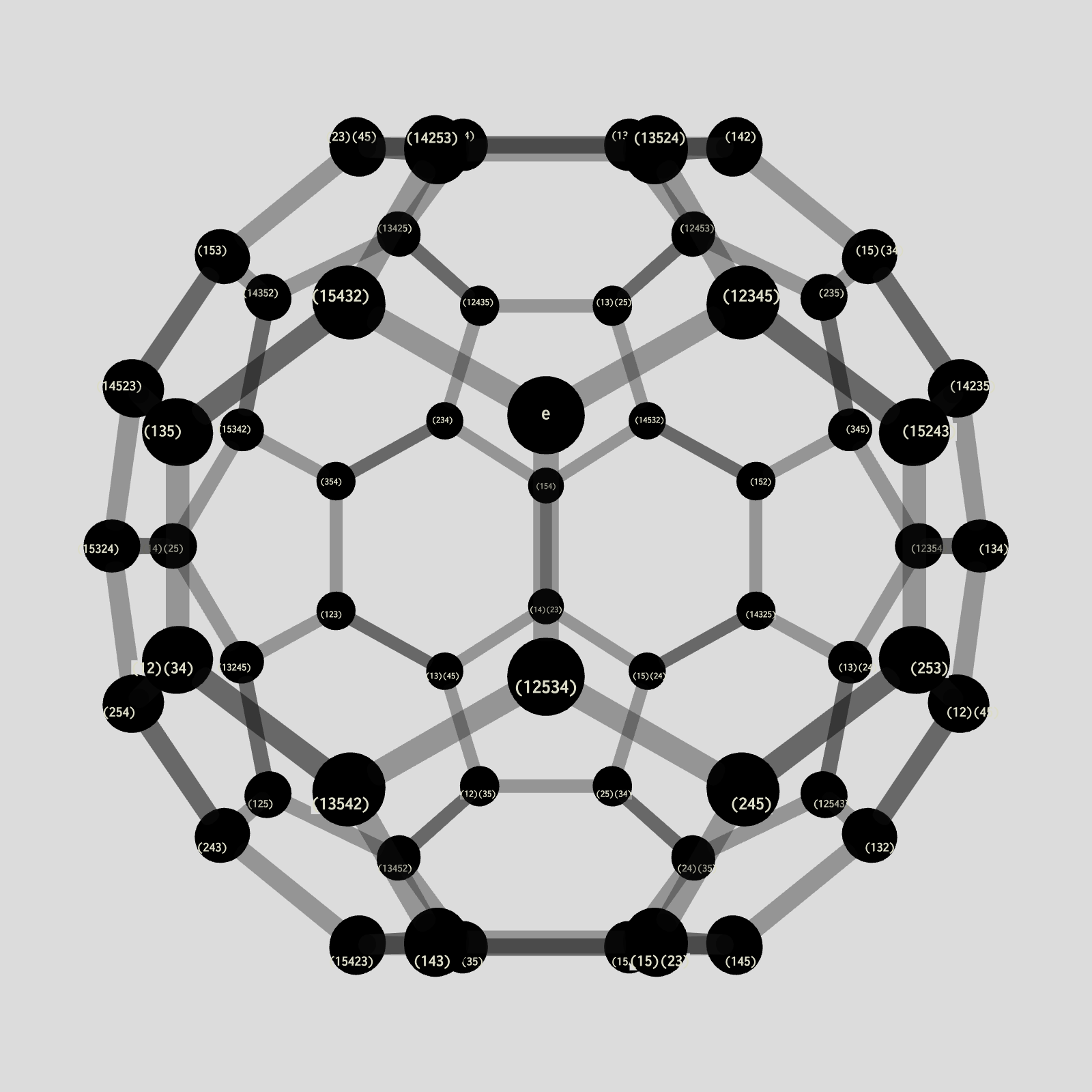}
    \caption{Labelling of carbons of Buckminsterfullerene with elements of $A_5$. Double bonds are not shown in the diagram, but there will be one double bond and two single bonds incident to each carbon. Scan the QR code for an interactive model.}
    \label{fig:bucky}
\end{figure}

A bijection is shown by the labelling of carbon atoms in Figure \ref{fig:bucky}. Let $v$ denote the vertex labelled $e$, and for each element $g \in A_5$ the corresponding vertex is $g(v)$.

\section{The Regular Representation}

We define the \emph{group algebra} $\CC[G]$ to be the algebra generated by elements $(e_g)_{g \in G}$ with the relations $e_g e_h = e_{gh}$ for each $g,h \in G$. Representations of $G$ extend linearly to representations of $\CC[G]$, and representations of $\CC[G]$ restrict to representations of $G$, so that representations of $G$ and $\CC[G]$ are equivalent. Let $M$ be the adjacency matrix of a graph $X=(V,E)$ with an f.t.a. action of $G$. Under the bijection of $V$ and $G$, Equation (\ref{eq:adjacency matrix def}) shows that $M$ acts on $\CC[G]$ by
\begin{equation}
	e_g \mapsto \sum_{(g,h) \in E} e_h
\end{equation}

Let $L_g$ and $R_g$ denote the left and right multiplication operators of $\CC[G]$ on itself by $e_g \in \CC[G]$. These two actions are called the left and right regular representations, respectively. Since $G$ acts by automorphisms, we have $L_k M = M L_k$, because for each $g \in G$

\begin{equation}
	L_k M e_g = \sum_{(g,h) \in E} e_{kh} = \sum_{(kg,kh) \in E} e_{kh} = M L_k e_g.
\end{equation}

Since $M e_{\id} = \sum_{(\id,h) \in E} L_h e_{\id}$, we have

\begin{align}
	M e_g &= M L_g e_{\id} = L_g M e_{\id} =  L_g \sum_{(\id,h) \in E} L_h e_{\id} \\
	&= \sum_{(\id,h) \in E} e_{gh} = e_g \sum_{(\id,h) \in E} R_h.
\end{align}

Therefore the left action of $M$ is equal to the right action of $\sum_{(\id,h) \in E} R_h$. More generally, for any $G$-invariant operator, $M'$, the left action of $M'$ is equal to right multiplication by the element $M' e_{\id}$, which we denote by $v(M')$. It is not hard to see that $v$ is an isomorphism $\End_G(\CC[G]) \to \CC[G]$, so you can show that for any $M' \in \End_G(\CC[G])$, $\sigma(M') = \sigma(v(M'))$. Thus we have reduced the problem of finding the eigenvalues of $H$ to finding those of
\begin{equation}
v(H) = v(I + M) = I+v(M) = I + \sum_{(\id,h) \in E} e_h.
\end{equation}

We have enough information to get a computer to calculate the matrix for $H$ using the regular representation, and this is done in \textbf{Notebook 1}. However, we can actually calculate the eigenvalues of $H$ without computing the matrix itself, using the characters of $G$. The decomposition of the regular representation tells us that

\begin{equation}
	\CC[G] \overset{\cong}{\longrightarrow} \bigoplus_{\lambda} V_\lambda^{d_\lambda}.
\end{equation}

where $\lambda$ indexes the irreducible representations $(\rho_\lambda,V_\lambda)$ of $G$, and the degree of $V_\lambda$ is $d_\lambda$. Let $\sigma(M)$ denote the multiset of eigenvalues of a matrix. For any left operator $M \in \End(\CC[G])$ you should be able to show that $\sigma(M) = \sigma(v(M))$. The decomposition of $\CC[G]$ into irreducible representations implies that

\begin{equation}
	\sigma(M) = \sum_\lambda d_\lambda\sigma(\rho_\lambda(v(M))),
\end{equation} 

where $\sum$ denotes the amalgamation of multisets. Now we can leverage the characters of $G$ to find the eigenvalues of $\rho_\lambda(v(M))$ for each irreducible representation $\rho_\lambda$. Let $\chi_\lambda$ denote the character of $V_\lambda$, i.e. $\chi_\lambda = \tr \rho_\lambda$. For each integer $k$ and any element $f \in \CC[G]$, we can calculate $\chi_\lambda(f^k) $ by expanding $f^k$ in terms of the basis vectors and using the additivity of the trace,

\begin{equation}
\chi_\lambda(f^k) = \chi_\lambda \left(\sum_{g \in G} c_g e_g \right) = \sum_{g \in G} c_g \chi_\lambda(g) .
\end{equation}

The only remaining step is to go from traces $\chi_\lambda(f^k)$ to the eigenvalues of $\rho_\lambda(f)$. 

\section{Newton's Formula}

In this section we determine the eigenvalues, $\sigma(M)$, for a matrix $M \in \Mat(\CC^d)$ given the values $\tr(M^k)$ for $k \in \{0, \ldots, d \}$. You probably know the following facts about the characteristic polynomial of $M$, denoted here by $f$.

\begin{itemize}
\item It can be written in terms of the eigenvalues $\sigma(M) = [\mu_1, \ldots, \mu_d]$ as

\begin{equation}
f(t) = \prod^d_{i=1} (t - \mu_i) = \sum_{i=0}^d (-1)^{d-i} e_{d-i}(\sigma(M)) t^i,
\end{equation}
where $e_i$ is the $i$\textsuperscript{th} elementary symmetric polynomial of $d$ variables, defined by
\begin{equation}
	e_i([x_1, \ldots, x_d ]) = \sum_{k_1 < \cdots < k_i} \prod^i_{j=1} x_{k_j}.
\end{equation}
We define $e_0 = 1$. The fact that polynomials $e_i$ are \emph{symmetric} justifies letting them take a multiset as input.
\item For any $k \in \NN$,
\begin{equation}
\tr(M^k) = \sum_{\mu \in \sigma(M)} \mu^k = p_k(\sigma(M)).
\end{equation}
The polynomials $p_k$ are called the power sum symmetric polynomials.

\end{itemize}

Now, the crucial fact is the following formula, which is very elegantly presented in \cite{macdonald}, and is due to Isaac Newton:

\begin{equation}
	n e_n = \sum_{r=1}^n (-1)^{r-1} p_r e_{n-r}.
\end{equation}

This allows us to compute the elementary symmetric polynomials $e_r$, $r \leq n$ recursively, using the power sum symmetric polynomials $p_r$, $r \leq n$. To recap: we compute $f(t)$ from the $e_n(\sigma(M))$, $n \leq d$; we compute $e_n(\sigma(M))$ from the $p_r(\sigma(M))$, $r \leq n$; and $p_r(\sigma(M))$ is just the trace $\tr(M^r)$. Since the roots of $f(t)$ are $\sigma(M)$, we are done. This is implemented in \textbf{Notebook 2}.

\section{Benzene and Other Aromatics}

In this section I will illustrate the method applied to benzene and other cyclic conjugated hydrocarbons, although elementary group theory suffices for the results of these examples. Consider the f.t.a. action of $\ZZ/6\ZZ$ on the benzene molecule generated by a $\pi/3$ radians rotation, denoted by $c$. We can extend this linearly to an action on $\CC^6$ by assigning a basis vector to each carbon. In matrix form, the action of $c$ is given by 

\begin{equation} \label{eq:Z/6Z gen matrix}
 c = \begin{bmatrix}
0 & 0 & 0 & 0 & 0 & 1 \\
1 & 0 & 0 & 0 & 0 & 0 \\
0 & 1 & 0 & 0 & 0 & 0 \\
0 & 0 & 1 & 0 & 0 & 0 \\
0 & 0 & 0 & 1 & 0 & 0 \\
0 & 0 & 0 & 0 & 1 & 0 
\end{bmatrix}.
\end{equation}

By choosing some clockwise labelling of the carbons $1,c,c^2,c^3,c^4,c^5$, it is clear that $\ZZ/6\ZZ$ acts on the carbons by group multiplication, so the action on $\CC^6$ is the regular representation of $\ZZ / 6\ZZ$. Since $\ZZ / 6\ZZ$ is Abelian, its irreducible representations are just the homomorphisms $\ZZ/6\ZZ \to \CC^*$. You should be able to work out that all such homomorphisms are of the form $\rho_{\lambda}: c \mapsto \zeta^{\lambda}$, where $\zeta = e^{\pi i/3}$ and $\lambda = 0, \ldots, 5$. Since the irreducible representations are one dimensional, $\chi_\lambda = \tr \rho_\lambda = \rho_\lambda$ so the character table is given directly by the previous statement, in Table \ref{tab:character table Z/6Z}.

\begin{table}[ht] 
\begin{tabular}{lllllll}
Size   		& $1$  	& $1$ 		& $1$   		& $1$  		& $1$   	& $1$        	\\
$\ZZ/6\ZZ$ 	& $1$ 	& $c$ 		& $c^2$ 		& $c^3$		& $c^4$ 	& $c^5$        	\\ \hline
$\chi_0$   	& $1$  	& $1$      	& $1$          	& $1$      	& $1$      	& $1$          	\\
$\chi_1$   	& $1$  	& $\zeta$	& $\zeta^2$  	& $\zeta^3$	& $\zeta^4$	& $\zeta^5$    	\\
$\chi_2$   	& $1$  	& $\zeta^2$	& $\zeta^4$  	& $1$		& $\zeta^2$	& $\zeta^4$    	\\
$\chi_3$   	& $1$  	& $\zeta^3$	& $1$		  	& $\zeta^3$	& $1$		& $\zeta^3$    	\\
$\chi_4$   	& $1$  	& $\zeta^4$	& $\zeta^2$  	& $1$		& $\zeta^4$	& $\zeta^2$    	\\
$\chi_5$	& $1$  	& $\zeta^5$	& $\zeta^4$  	& $\zeta^3$	& $\zeta^2$	& $\zeta^1$    	\\
\end{tabular}
\caption{The character table of $\ZZ/6\ZZ$. The rows index irreducible representation and the columns index conjugacy classes.}
\label{tab:character table Z/6Z}
\end{table} 

The adjacency matrix of $M$ is written in Equation (\ref{eq:benzene huckel matrix}). By comparison with Equation (\ref{eq:Z/6Z gen matrix}) we conclude that $M = c + c^{-1}$. Since the regular representation decomposes into a direct sum of the irreducible representations $\rho_\lambda$, the eigenvalues of $M$ are simply the union of the eigenvalues of $\rho_{\lambda}(M)$. Since $\rho_{\lambda}(M)$ are scalars, the eigenvalues are simply 
\begin{align}
	\mu_\lambda &= \rho_{\lambda}(M) \\
	&= \rho_{\lambda}(c) + \rho_{\lambda}(c^{-1}) = \chi_\lambda(c) + \chi_\lambda(c^{-1}) \\
	&= \zeta^\lambda + \overline{\zeta}^\lambda \\
	&= 2 \cos(\pi \lambda / 3).
\end{align}
Therefore the eigenvalues with multiplicity are $-2,-1,-1,1,1,2$. Therefore the eigenvalues of the H\"uckel matrix $H = \alpha I + \beta M$ are
\begin{equation}
	\alpha -2 \beta, \alpha - \beta, \alpha - \beta, \alpha + \beta, \alpha + \beta, \alpha + 2\beta.
\end{equation}

Setting $\alpha = 0$ and $\beta = -1$, the ground state energy is $2(-2) + 2(-1) + 2(-1) = -8$. The same analysis can be applied to any cyclic conjugated hydrocarbon, and a general formula for the ground state energy of the molecule of length $n \in \ZZ_{\geq 3}$ can be obtained,
\begin{equation} \label{eq:4n+2 rule}
E_n = \begin{cases}
	-4 \cos(\pi/n) / \sin(\pi/n) , &\text{ if } n \equiv 0 \pmod{4} \\
	-4/\sin(\pi/n ), &\text{ if } n \equiv 2 \pmod{4}
\end{cases}.
\end{equation}
Note that $n$ must be even because of the alternating single and double bonds, and $n$ must be at least three. Since $\cos(\pi/n) \in (0,1)$ for all $n \in \NN$, Equation (\ref{eq:4n+2 rule}) implies that the ground state energy of cyclic conjugated hydrocarbons of length $n \equiv 2 \pmod{4}$ is lower than $n \equiv 0 \pmod{4}$, and therefore they are more likely to form. This is called the ``$4p+2$ rule'' in chemistry.

\section{The buckyball spectrum}

In this final section, we discuss the previous result in the specific context of Buckminsterfullerene. As mentioned, the Buckminsterfullerene graph has a f.t.a. action of the alternating group $G=A_5$. Figure \ref{fig:bucky} shows that the single bond adjacency matrix of Buckminsterfullerene is represented by the element $M_1 = (12345) + (15432)$ and the double bond adjacency matrix is represented by $M_2 = (12)(34)$. Thus the Hückel Matrix is of the form

\begin{equation}
H=\alpha I + \beta M_1 + \gamma M_2 = \alpha \id + \beta ((12345) + (15432)) + \gamma (12)(34).
\end{equation}

From this equation, we calculate $H^k$ for $1 \leq k \leq 5$ as an element of $\CC[G]$.

The character table of $A_5$ is shown in Table \ref{tab:A5 char table}, which allows us to determine $\chi_\lambda(H^k)$. For example, the representation $Y$ has degree three, so we first calculate the expansions of $H^0, \ldots H^3$ in $\CC[G]$ using the group law of $G$. For $\alpha =0, \beta =-1, \gamma =-1$, these expansions come out as
\begin{align*}
	H^0 &= \id, \quad H^1 = (12)(34) + (12345) + (15432)\\
	H^2 &= 3\id + (245) + (254) + (135) + (13524) + (14253) + (153) \\
	H^3 &= (23)(45) + 5(12)(34) + 5(12345) + 2(12534) + (13524) + (13254) \\
	&\phantom{=}+ (14352) + (14523) + (14235)+ (14253) + (15432) + (15)(23) + (15324). \\
\end{align*}
Now we use the character table to compute $\chi_Y(H^k)$. For example,
\begin{equation}
	\chi_Y(H^2) = 3(3) + 0 + 0 + 0 + \frac{1-\sqrt{5}}{2} + \frac{1-\sqrt{5}}{2} + 0 = 10 - \sqrt{5}.
\end{equation}

\begin{table}[ht]
\begin{tabular}{llllll}
Size      & $1$   & $20$        & $15$           & $12$                     & $12$                     \\
$A_5$ & $1$ & $x=(123)$ & $y=(12)(34)$ & $z=(12345)$            & $z^2=(13524)$          \\ \hline
$U$     & $1$   & $1$         & $1$            & $1$                      & $1$                      \\
$V$     & $4$   & $1$         & $0$            & $-1$                     & $-1$                     \\
$W$     & $5$   & $-1$        & $1$            & $0$                      & $0$                      \\
$Y$     & $3$   & $0$         & $-1$           & $\frac{1+\sqrt{5}}{2}$ & $\frac{1-\sqrt{5}}{2}$ \\
$Z$     & $3$   & $0$         & $-1$           & $\frac{1-\sqrt{5}}{2}$ & $\frac{1+\sqrt{5}}{2}$
\end{tabular} 
\caption{The character table of $A_5$. The rows index irreducible representation and the columns index conjugacy classes.}
\label{tab:A5 char table}
\end{table} 

We then calculate the elementary symmetric polynomials in terms of the power polynomials, and substitute in our values for the traces. This gives us the coefficients of the characteristic polynomial for each $\lambda$. For all but one of the representations, the roots of the characteristic polynomial of $\rho_{\lambda}(H)$ can be solved exactly since the degree is less than $5$, and in fact the degree $5$ polynomials factors, so that is also solvable exactly. This is implemented in \textbf{Notebook 2}. The characteristic polynomials for each irreducible representation, with the coefficients $\alpha = 0$, $\beta =1$, and $\gamma = r$, are as follows

\begin{align*}
f_U(t)&= r - t + 2 \\
f_Y(t)&= r^3 + r^2 t - r t^2 - t^3 - \sqrt{5} r t - \sqrt{5} t^2 + r t + t^2 + \frac{1}{2} \sqrt{5} r + \frac{3}{2} \sqrt{5} t - \frac{1}{2} r + \frac{1}{2} t + \sqrt{5} + 3 \\
f_Z(t)&= r^3 + r^2 t - r t^2 - t^3 + \sqrt{5} r t + \sqrt{5} t^2 + r t + t^2 - \frac{1}{2} \sqrt{5} r - \frac{3}{2} \sqrt{5} t - \frac{1}{2} r + \frac{1}{2} t - \sqrt{5} + 3 \\
f_V(t)&= (r^2 - t^2 + 2 r - t + 1) (r^2 - t^2 - t + 1) \\
f_W(t)&= (r^3 - r^2 t - r t^2 + t^3 - r^2 + 2 r t - t^2 + r - 3 t + 2) (r^2 - t^2 - r - t + 1). \\
\end{align*}

We can calculate the eigenvalues by taking the roots of these polynomials. For $\alpha=0, \beta=-1,\gamma=-1$, the resulting ground state energy is $-93.16$, and this is confirmed by the explicit matrix calculation.

\section*{Further Applications}

In this article, we used representation theory to find the $\pi$-electron energy levels of molecules with either cyclic or icosahedral symmetry. There are other finite subgroups of $SO(3)$ which may have an f.t.a. action on a conjugated molecule, but since conjugation requires a degree of local planarity, the only likely candidates are the larger ones (note that the icosahedral group is the largest non-cyclic example).

\begin{figure}
    \centering
    \includegraphics[width=0.5\linewidth]{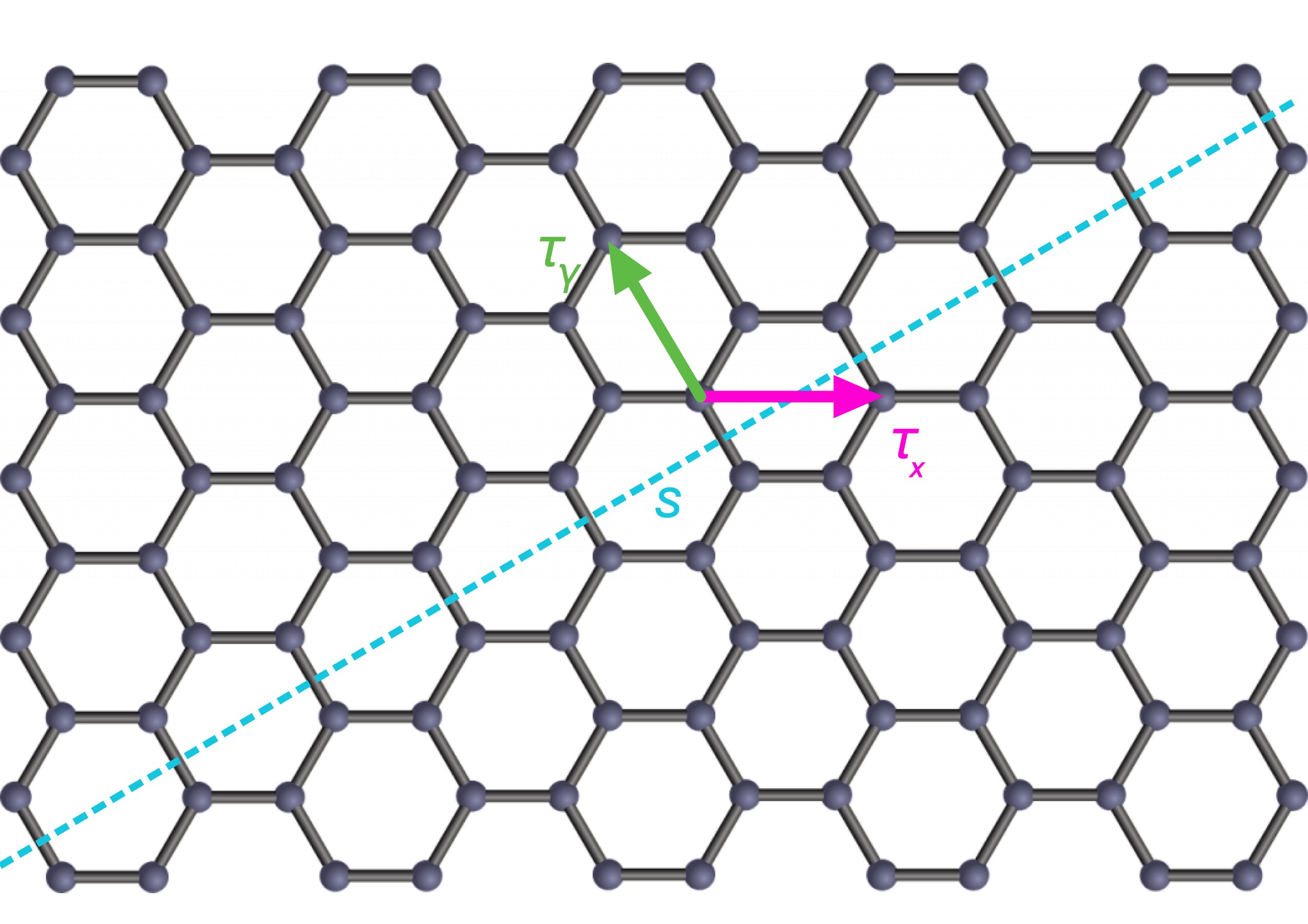}
    \caption{A diagram of the of the structure of a section of graphene, ignoring double bonds. The generators of the group action are shown: the arrows indicate translations and the dashed line indicates reflection. Adapted from \cite{graphene_image}.}
    \label{fig:graphene}
\end{figure}

More examples are available if we relax some of our conditions, such as the finiteness of the group $G$. The molecule graphene has a honeycomb skeleton extending in two dimensions. If we approximate graphene by an infinite honeycomb, it has an f.t.a. action by the group $G = \ZZ^2 \rtimes (\ZZ/2\ZZ)$, where $\ZZ^2$ is generated by the translations $\tau_x, \tau_y$ indicated in Figure \ref{fig:graphene} and $\ZZ / 2\ZZ$ is generated by the reflection $s$ along the dashed line. The adjacency matrix $M$ is such that $v(M) = s(1 + \tau_x + \tau_x^{-1} \tau_y^{-1})$. Instead of a finite decomposition, we have a decomposition into a continuous space of two-dimensional representations parametrised by the torus, and thus the spectrum is continuous, $\sigma(M) = [-3,3]$. I believe that the energy levels of de-localised electrons in three-dimensional metal lattices such as tantalum or antimony could be computed in the same way, even though they arise for chemical reasons quite different from conjugation. 

\section{Conclusion}
Finding the spectrum of H\"uckel matrices using representation theory has three advantages over the traditional method. Firstly, it averts the need to write out a large matrix by hand. Secondly, although a computer can quickly diagonalise a large finite matrix such as the $60 \times 60$ buckyball adjacency matrix, this brute force approach fails to explain why the eigenvalues have such high multiplicity, or why they are algebraic, and so on. Additionally our technique can be done fully by hand, and the parameters can be adjusted more easily. Thirdly, in the case of infinite molecules—or rather the limit of finite molecules which can extend without bound—direct computation is not possible.

\section*{Further Reading}
For a reference on the mathematics of quantum mechanics I recommend \cite{neumann}. For H\"uckel Theory, I recommend \cite{Coulson_Mallion_O’Leary_1978}. For representation theory, I recommend \cite{serre}. 

I would like to thank my twin brother, Oscar, for his help in clarifying the chemical aspects of this article.

\bibliographystyle{abbrv}
\bibliography{biblio}

\end{document}